\documentclass[12pt]{article}

\usepackage[T2A]{fontenc}
\usepackage[cp1251]{inputenc}
\usepackage{amsfonts}
\usepackage{amssymb}
\usepackage{amsmath}

\usepackage[usenames]{color}

\begin{document}
	
	\newtheorem{theorem}{Theorem}
	\newtheorem{lemma}{Lemma}
	\newtheorem{corollary}{Corollary}
	\newtheorem{definition}{Definition}
	\newtheorem{proposition}{Proposition}
	\newtheorem{remark}{Remark}

\begin{center}

\textbf{\Large The Klebanov theorem \\ for the group $\mathbb{R}\times \mathbb{Z}(2)$}

\bigskip

{\Large Margaryta Myronyuk}

\bigskip

\textit{B. Verkin Institute for Low Temperature Physics and Engineering\\ of the National Academy of
	Sciences of Ukraine, \\ Nauky Ave. 47, Kharkiv, Ukraine}

\medskip

myronyuk@ilt.kharkov.ua

\medskip

ORCID 0000-0001-9596-9292
\end{center}

\begin{abstract}
L. Klebanov proved the following theorem. Let $\xi_1, \dots, \xi_n$ be independent random variables. Consider linear forms $L_1=a_1\xi_1+\cdots+a_n\xi_n,$ 
$L_2=b_1\xi_1+\cdots+b_n\xi_n,$
$L_3=c_1\xi_1+\cdots+c_n\xi_n,$
$L_4=d_1\xi_1+\cdots+d_n\xi_n,$  
where the coefficients $a_j, b_j, c_j, d_j$ are real numbers. If the random vectors $(L_1,L_2)$ and $(L_3,L_4)$ are identically distributed, then all $\xi_i$ for which $a_id_j-b_ic_j\neq 0$ for all $j=\overline{1,n}$ are Gaussian random variables. The present article is devoted to an
analogue of the Klebanov theorem in the case when random variables take
values in the group $\mathbb{R}\times \mathbb{Z}(2)$ and the coefficients of
the linear forms are topological endomorphisms of this group.
\end{abstract}



\bigskip
	
	\section{Introduction}\label{sec1}
	
	\medskip
	
Characterization theorems in mathematical statistics provide conditions under which the distributions of random variables can be uniquely determined from specified properties of functions of these variables. Such theorems on the real line  are well studied (\cite{KaLiRa}). In particular, Darmois (\cite{Darmois}) and Skitovich (\cite{Skitovich}) established that the independence of two linear forms of 
$n$ independent random variables characterizes the Gaussian distribution on the real line. Heyde (\cite{Heyde}) obtained a similar results: the symmetry of the conditional distribution of one linear form given another also characterizes the Gaussian distribution. Klebanov (\cite{Klebanov}) proved the following theorem.

\medskip
	
\noindent \textbf{The Klebanov theorem.} \textit{ Let $\xi_1, \dots, \xi_n$ be independent random variables. Consider linear forms $L_1=a_1\xi_1+\cdots+a_n\xi_n,$ 
		$L_2=b_1\xi_1+\cdots+b_n\xi_n,$
		$L_3=c_1\xi_1+\cdots+c_n\xi_n,$
		$L_4=d_1\xi_1+\cdots+d_n\xi_n,$  
		where the coefficients $a_j, b_j, c_j, d_j$ are real numbers. If the random vectors $(L_1,L_2)$ and $(L_3,L_4)$ are identically distributed, then all $\xi_i$ for which $a_id_j-b_ic_j\neq 0$ for all $j=\overline{1,n}$ are Gaussian random variables. }
	
\medskip

\noindent The Klebanov theorem implies both the Darmois-Skitovich theorem and the Heyde theorem.

A large number of papers are devoted to the study of linear forms of independent random variables with values in locally compact Abelian groups. The main attention is given to analogues of the Darmois–Skitovich and Heyde theorems. To the best of our knowledge, the Klebanov theorem was generalized to groups only in \cite{My2024}, where the coefficients of the linear forms are integers.

The purpose of the present paper is to establish an analogue of the Klebanov theorem in the case where the random variables take values in the locally compact Abelian group 
$X=\mathbb{R}\times\mathbb{Z}(2)$ ($\mathbb{Z}(2)$ is the cyclic group of order 2) and the coefficients of the linear forms are continuous endomorphisms. The group $X=\mathbb{R}\times\mathbb{Z}(2)$
is important because it is the simplest nontrivial example of a non-connected locally compact Abelian group whose structure simultaneously involves a connected component and a discrete component. This combination makes $X$ a testing ground for many problems in harmonic analysis, probability theory on groups, and characterization theorems (\cite{Fe2022}, \cite{Fe2024}, \cite{Trukhina}, \cite{Zolotarev}).

It should be noted that a continuous endomorphism of the real line is given by multiplication by a real number; if this number is nonzero, the endomorphism is an automorphism. In contrast, for arbitrary locally compact Abelian groups the situation is substantially more intricate: the set of continuous endomorphisms may be infinite, whereas the subset of automorphisms may be finite. Consequently, generalizations of the Klebanov theorem to arbitrary groups present significant difficulties, but they may also give rise to distributions of a rather unexpected nature.

\section{Notation and definitions}\label{sec2}
	
The article employs standard results from abstract harmonic analysis
	(see e.g. \cite{HeRo1}).
	
Let $X$ be a second countable locally compact  Abelian group,$Y=X^\ast$ be its character group, and  $(x,y)$ be the value of a
character $y \in Y$ at an element $x \in X$. Denote by $Aut(X)$ the group of all topological automorphisms of the group $X$. Denote by $End(X)$ the ring of all continuous endomorphisms of $X$. Let $H$ be a subgroup
of $Y$. Denote by $A(X,H)=\{x \in X: (x,y)=1 \ \ \forall \ y \in
H\}$ the annihilator of $H$. Let $\delta\in End(X)$. The adjoint endomorphism $\widetilde{\delta}$ is defined by the formula $(\delta x,y)=(x, \widetilde{\delta}y)$ for all $x\in X$, $y\in Y$. 
	
	Let ${M^1}(X)$ be the convolution semigroup of probability
	distributions on $X$,  $$\widehat \mu(y) = \int_X (x, y) d\mu(x)$$ be
	the characteristic function of a distribution $\mu \in {M^1}(X)$,
	and $\sigma(\mu)$ be the support of $\mu$. If $H$ is a closed
	subgroup of $Y$ and $\widehat \mu(y)=1$ for $y \in H$, then
	$\widehat\mu(y+h) = \widehat\mu(y)$ for all $y \in Y$, $ h \in H$ and
	$\sigma(\mu) \subset A(X, H)$. For $\mu \in {M^1}(X)$ we define the
	distribution $\bar \mu \in M^1(X)$ by the rule $\bar \mu(B) =
	\mu(-B)$ for all Borel sets $B \subset X$. Note that $\widehat
	{\bar \mu}(y) = \overline{\widehat \mu(y)}$.
	
	Let $x\in X$. Denote by $E_x$ the degenerate distribution
	concentrated at the point $x$, and by $D(X)$ the set of all
	degenerate distributions on
	$X$. A distribution $\gamma \in {M^1}(X)$ is called Gaussian (\cite[Ch.
	4.6]{Pa}) if its characteristic function can be represented in the
	form
	\begin{equation*}\label{i1}
		\widehat\gamma(y)= (x,y)\exp\{-\varphi(y)\},
	\end{equation*}	
	where $x \in X$ and $\varphi(y)$ is a continuous nonnegative
	function satisfying the equation
	\begin{equation*}\label{i2}
		\varphi(u+v)+\varphi(u-v)=2[\varphi(u)+
		\varphi(v)], \quad u, \ v \in Y.
	\end{equation*}
	Denote by $\Gamma(X)$ the set
	of Gaussian distributions on $X$. We note that according to this
	definition $D(X)\subset \Gamma(X)$. Denote by $I(X)$ the set of
	shifts of Haar distributions $m_K$ of compact subgroups $K$ of the
	group $X$. Note that
	
	\begin{equation}
		\widehat{m}_K(y)=
		\left\{%
		\begin{array}{ll}
			1, & \hbox{$y\in A(Y,K)$;} \\
			0, & \hbox{$y\not\in A(Y,K)$.} \\
		\end{array}%
		\right.
	\end{equation}
	We note that if a distribution $\mu \in \Gamma(X)*I(X)$, i.e.
	$\mu=\gamma*m_K$, where $\gamma\in \Gamma(X)$, then $\mu$ is
	invariant with respect to the compact subgroup $K \subset X$ and
	under the natural homomorphism $X \mapsto X/K \ \mu$ induces a
	Gaussian distribution on the factor group $X/K$. Therefore the class
	$\Gamma(X)*I(X)$ can be considered as a natural analogue of the
	class $\Gamma(X)$ on locally compact Abelian groups. 
	
	Let $f(y)$ be a function on $Y$, and $h\in Y.$ Denote by $\Delta_h$
	the finite difference operator
	$$\Delta_h f(y)=f(y+h)-f(y).$$
	A function $f(y)$ on $Y$ is called a polynomial if
	$$\Delta_{h}^{l+1}f(y)=0$$
	for some $l$ and for all $y,h \in Y$. The minimal $l$ for which  this equality holds is called the degree of the polynomial $f(y)$.
	
	An integer $a$ is said to be admissible for a group $X$ if the set $\{ax, x\in X\}
	\ne \{0\}$. The admissibility of integers $a$
	can be considered as a group analogue of the condition  $a\ne 0$ for the case of  $X = \mathbb{R}$.
	
\section{Main results}\label{sec3}

Let $X=\mathbb{R}\times\mathbb{Z}(2)$. Denote by $x = (t, k)$, $t \in \mathbb{R}$, $k \in \mathbb{Z}(2)$, elements of the group $X$. Let $Y$ be the character group of the group $X$. The group $Y$ is topologically isomorphic to the group $\mathbb{R}\times \mathbb{Z}(2)$.
Denote by $y = (s, l)$, $s \in \mathbb{R}$, $l \in \mathbb{Z}(2)$, elements of the group $Y$. Every continuous endomorphism
$a$ is of the form $a(t, k) = (a' t, a'' k)$, where $a'\in \mathbb{R}$, $a''\in \{0,1\}$. If $a\in Aut(X)$ then $a' \neq 0$, $a''=1$. It is obvious that $a = \widetilde{a}$.

\begin{definition}\label{ClassTheta}
	Let $X=\mathbb{R}\times \mathbb{Z}(2)$ and $\mu\in M^1(X)$. We say that $\mu\in\Theta(X)$ if the characteristic function $\widehat{\mu}(s,l)$ is represented in the form 
	\begin{equation}\label{def1}
		\widehat\mu(s,l)=
		\left\{%
		\begin{array}{ll}
			\exp\{-\sigma s^2+ i\beta s\}, & \hbox{$s\in \mathbb{R}, l=0$;} \\
			\kappa \exp\{-\sigma' s^2+ i\beta' s\}, & \hbox{$s\in \mathbb{R}, l=1$,} \\
		\end{array}%
		\right.
	\end{equation}
	where either $0 < \sigma' < \sigma$, and $0< |\kappa| \leq \sqrt{\sigma'\over \sigma} \exp\left\lbrace -{(\beta-\beta')^2\over 4(\sigma-\sigma')} \right\rbrace $ or $\sigma=\sigma'\geq 0$,  $\beta=\beta'$, $ |\kappa| \leq 1$.
\end{definition}
Note that the class of distributions $\Theta(X)$ first appeared in \cite{Fe2020}.

\begin{definition}
	Let $X=\mathbb{R}\times \mathbb{Z}(2)$ and $\mu\in M^1(X)$. We say that $\mu\in \Lambda(X)$ if $\widehat{\mu}(s,0)$ is the characteristic function of a Gaussian distribution on
	the real line, possibly degenerate.
\end{definition}

Let $X$ be a second countable locally compact Abelian group. Let $\xi_1, \dots, \xi_n$ be independent random variables with values in $X$. Consider the linear forms 

\begin{equation}\label{th1.1.1}
	L_1=a_1\xi_1+\cdots+a_n\xi_n,
\end{equation}
\begin{equation}\label{th1.1.2}
	L_2=b_1\xi_1+\cdots+b_n\xi_n,
\end{equation}
\begin{equation}\label{th1.1.3}
	L_3=c_1\xi_1+\cdots+c_n\xi_n,
\end{equation}
\begin{equation}\label{th1.1.4}
	L_4=d_1\xi_1+\cdots+d_n\xi_n,
\end{equation}
where the coefficients $a_j, b_j, c_j, d_j\in End(X)$. We will keep the designation $L_1$, $L_2$, $L_3$ and $L_4$ throughout the article. 

\begin{theorem}\label{main.nonzero}
	Let $X=\mathbb{R}\times \mathbb{Z}(2)$.	
	Let $\xi_1, \dots, \xi_n$ be independent random variables with values in $X$ and distributions $\mu_{\xi_j}$ with non-vanishing characteristic functions. Consider the linear forms $(\ref{th1.1.1})-(\ref{th1.1.4})$, 
	where $a_j, b_j, c_j, d_j\in End(X)$. Suppose that the random vectors $(L_1,L_2)$ and $(L_3,L_4)$ are identically distributed. We have the following possibilities:
	
	$(1)$ for all those $i$ for which $a_i d_j - b_i c_j\in Aut(X)$ for all  $j=\overline{1,n}$ distributions $\mu_{\xi_i}\in\Gamma(X)$;
	
	$(2)$ for all those $i$ for which $a'_i d'_j - b'_i c'_j\neq 0$ for all  $j=\overline{1,n}$ and at least one of the numbers $a''_i$, $b''_i$ is equal to 1 distributions $\mu_{\xi_i}\in\Theta(X)$;
	
	$(3)$ for all those $i$ for which $a'_i d'_j - b'_i c'_j\neq 0$ for all  $j=\overline{1,n}$ distributions $\mu_{\xi_i}\in\Lambda(X)$.
	
\end{theorem}

To prove Theorem \ref{main.nonzero} we need some lemmas. The following lemma was proved earlier in \cite{My2024} in the case when the coefficients of the linear forms $(\ref{th1.1.1})-(\ref{th1.1.4})$ are integers. The proof follows the same structure and is provided here for completeness.

\begin{lemma}\label{Equation}
	Let $X$ be a second countable locally
	compact Abelian group, $Y=X^*$. Let $\xi_1, \dots, \xi_n$ be independent random variables with values in $X$ and distributions $\mu_{\xi_j}$. Consider the linear forms $(\ref{th1.1.1})-(\ref{th1.1.4})$. The random vectors $(L_1,L_2)$ and $(L_3,L_4)$ are identically distributed if and only if the characteristic functions $\widehat\mu_{\xi_j}(y)$ satisfy the equation
	\begin{equation}\label{Eq1}
		\prod_{j=1}^{n}  \widehat\mu_{\xi_j} (\widetilde{a}_j u+ \widetilde{b}_j v)=\prod_{j=1}^{n}  \widehat\mu_{\xi_j} (\widetilde{c}_j u+ \widetilde{d}_j v),\quad u,v\in Y.
	\end{equation}
\end{lemma}

\textbf{Proof.} Taking into account the independence of $\xi_1, ..., \xi_n$, we obtain from the definition that
\begin{eqnarray*}
	\widehat\mu_{(L_1,L_2)}(u,v)=\mathbf{E}\left[ (L_1,L_2) (u,v) \right] =\mathbf{E}\left[ (L_1,u) (L_2,v) \right] = \\ =\mathbf{E}\left[ (a_1\xi_1+\cdots+a_n\xi_n,u) (b_1\xi_1+\cdots+b_n\xi_n,v) \right] = \\ =\mathbf{E}\left[ (a_1\xi_1,u)\cdots(a_n\xi_n,u) (b_1\xi_1,v)\cdots(b_n\xi_n,v) \right] = \\ =\mathbf{E}\left[ (\xi_1,\widetilde{a}_1 u)\cdots(\xi_n, \widetilde{a}_n u) (\xi_1,\widetilde{b}_1 v)\cdots(\xi_n, \widetilde{b}_n v) \right]= \\ =\mathbf{E}\left[ (\xi_1,\widetilde{a}_1 u +\widetilde{b}_1 v)\cdots(\xi_n, \widetilde{a}_n u+\widetilde{b}_n v) \right]= \\=\prod_{j=1}^{n} \mathbf{E}\left[ (\xi_j,\widetilde{a}_j u +\widetilde{b}_j v) \right]  = \prod_{j=1}^{n}  \widehat\mu_{\xi_j} (\widetilde{a}_j u+ \widetilde{b}_j v),\quad u,v\in Y.
\end{eqnarray*}
In the same way we obtain that
\begin{equation*}
	\widehat\mu_{(L_3,L_4)}(u,v)=\prod_{j=1}^{n}  \widehat\mu_{\xi_j} (\widetilde{c}_j u+ \widetilde{d}_j v),\quad u,v\in Y.
\end{equation*} 
The vectors $(L_1,L_2)$ and $(L_3,L_4)$ are identically distributed if and only if $\widehat\mu_{(L_1,L_2)}(u,v)=\widehat\mu_{(L_3,L_4)}(u,v)$. Thus, equation
(\ref{Eq1}) is valid.
\hfill $\square$

\bigskip

We need the following group analogue of the Cramer decomposition theorem for Gaussian distributions.

\begin{lemma}\label{GrKramer} $(\cite[\S5.22]{FeBook1})$ Let $X$ be a second countable locally compact Abelian group containing no subgroup topologically isomorphic to the circle group. Let $\gamma\in \Gamma(X)$ and $\gamma=\gamma_1*\gamma_2$, where $\gamma_1, \gamma_2 \in M^1(X)$. Then $\gamma_1, \gamma_2 \in \Gamma(X)$.
\end{lemma}

The following lemma is the group analogue of the Marcinkiewicz theorem.

\begin{lemma}\label{GrMarc} $(\cite[\S5.11]{FeBook2})$
	Let $X$ be a second countable locally compact Abelian group containing no subgroup topologically isomorphic to the circle group. Let $\mu \in M^1(X)$ and the characteristic function $\widehat{\mu}(y)$ is of the form
	\begin{equation*}
		\widehat{\mu}(y)=\exp\{ \varphi(y)\}, \quad y\in Y,
	\end{equation*}
	where $\varphi(y)$ is a continuous polynomial. Then $\mu\in \Gamma(X)$.
\end{lemma}

\begin{lemma}\label{AlmostKramer} $(\cite{Fe2020})$ Suppose that $\mu\in \Theta(\mathbb{R}\times \mathbb{Z}(2))$ and the characteristic function of $\mu$ does not vanish. Let $\mu=\mu_1*\mu_2$, where $\mu_1, \mu_2 \in M^1(\mathbb{R}\times \mathbb{Z}(2))$. Then $\mu_1, \mu_2 \in \Theta(\mathbb{R}\times \mathbb{Z}(2))$.
\end{lemma}

Lemma \ref{AlmostKramer} can be considered as an analogue of the Cramer decomposition theorem for distributions which belong to $\Theta(\mathbb{R}\times \mathbb{Z}(2))$.

\begin{lemma}\label{ClassTheta} $(\cite{Fe2020})$
	Let $X=\mathbb{R}\times \mathbb{Z}(2)$, and let $f(s,l)$ be a function on the group $Y$ of the form (\ref{def1}),
	where $\sigma, \sigma' \geq 0$, $\beta, \beta',\kappa \in \mathbb{R}$. Then $f(s,l)$ is the characteristic function of a signed measure $\mu$ on the group $X$. Moreover, $\mu\in M^1(X)$ if and only if either $0 < \sigma' < \sigma$ and $0< |\kappa| \leq \sqrt{\sigma'\over \sigma} \exp\left\lbrace -{(\beta-\beta')^2\over 4(\sigma-\sigma')} \right\rbrace $ or $\sigma=\sigma'$,  $\beta=\beta'$, $ |\kappa| \leq 1$. In the last case $\mu\in\Gamma(\mathbb{R})*M^1(\mathbb{Z}(2))$.
\end{lemma}

\begin{lemma}\label{EntireFunction} (\cite[Lemma 6.9]{FeBook1}) 
	Let $X = \mathbb{R}\times G$, where $G$ is a locally compact Abelian group. Let $Y=X^*$ and $H=G^*$. Denote by $(s, h)$, $s \in \mathbb{R}$, $h \in H$,
	elements of the group $Y$. Let $\mu\in M^1(X)$ and assume that the characteristic function $\widehat{\mu}(s, 0)$
	can be extended to the complex plane as an entire function in $s$. Then for each fixed $h\in H$ the function $\widehat{\mu}(s, h)$ can be also extended to the complex plane as an entire function in $s$.
	Moreover, the following inequality
	\begin{equation}\label{EntireFunctionIneq}
		\max_{|s|\leq r} |\widehat\mu(s,h)| \leq \max_{|s|\leq r} |\widehat\mu(s,0)|, \quad h\in H,
	\end{equation}
	holds.
\end{lemma}

The following proposition is a special case of Theorem 1 of \cite{My2024}.

\begin{proposition}\label{z2} (\cite{My2024}).
	Let $X=\mathbb{Z}(2)$. Let $\xi_1, \dots, \xi_n$ be independent random variables with values in $X$ and distributions $\mu_{\xi_j}$ with non-vanishing characteristic functions. Consider the linear forms $(\ref{th1.1.1})-(\ref{th1.1.4})$, 
	where the coefficients $a_j, b_j, c_j, d_j\in \mathbb{Z}(2)$. If the random vectors $(L_1,L_2)$ and $(L_3,L_4)$ are identically distributed, then for all those $i$ for which
	\begin{equation}\label{th1.1}
		a_i d_j - b_i c_j =1 \text{ for all } j=\overline{1,n},
	\end{equation}
	$\mu_{\xi_i}\in D(X)$.	
\end{proposition}


\textit{Proof of Theorem \ref{main.nonzero}.} Since $X=\mathbb{R}\times \mathbb{Z}(2)$, we have $Y\approx\mathbb{R}\times \mathbb{Z}(2)$. Therefore $\widetilde{a}_j=a_j$, $\widetilde{b}_j=b_j$, $\widetilde{c}_j=c_j$, $\widetilde{d}_j=d_j$. By Lemma \ref{Equation} the characteristic functions $\widehat\mu_{\xi_j}(y)$ satisfy equation (\ref{Eq1}) which can be written in the following form

\begin{eqnarray}\label{Eq1-1}
	\prod_{j=1}^{n}  \widehat\mu_{\xi_j} (a'_j s_1+ b'_j s_2, a''_j l_1+ b''_j l_2 )
	=\prod_{j=1}^{n}  \widehat\mu_{\xi_j} (c'_j s_1+ d'_j s_2, c''_j l_1+ d''_j l_2 ),\nonumber \\ s_1,s_2\in\mathbb{R}, l_1, l_2 \in \mathbb{Z}(2).
\end{eqnarray}
Substituting $l_1=l_2=0$ in (\ref{Eq1-1}), we obtain
\begin{eqnarray}\label{Eq1-2}
	\prod_{j=1}^{n}  \widehat\mu_{\xi_j} (a'_j s_1+ b'_j s_2, 0)=\prod_{j=1}^{n}  \widehat\mu_{\xi_j} (c'_j s_1+ d'_j s_2, 0 ),\quad s_1,s_2\in\mathbb{R}.
\end{eqnarray}
It follows now from Lemma \ref{Equation} and the Klebanov theorem that 
\begin{eqnarray}\label{Eq1-2}
	\widehat\mu_{\xi_i} (s, 0)= e^{-\sigma_i s^2+ i \beta_i s} ,\quad s\in\mathbb{R},
\end{eqnarray}
where $\sigma_i\geq 0$, $\beta_i\in\mathbb{R}$, for all those $i$ for which $a'_i d'_j - b'_i c'_j\neq 0$ for all  $j=\overline{1,n}$. Thus, we proved statement (3) of Theorem \ref{main.nonzero}.

Now we prove statement (2) of Theorem \ref{main.nonzero}. Put $\nu_{\xi_j}=\mu_{\xi_j}*\bar{\mu}_{\xi_j}$. Hence $\widehat\nu_{\xi_j}(y)=|\widehat\mu_{\xi_j}(y)|^2 > 0$, $y\in Y$. Obviously, the characteristic functions $\widehat\nu_{\xi_j}(y)$ also satisfy equation (\ref{Eq1}). If we prove that either $\nu_{\xi_j}\in\Gamma(X)$ or $\nu_{\xi_j}\in\Theta(X)$, then applying either Lemma \ref{GrKramer} or Lemma \ref{AlmostKramer} respectively, we get that either $\mu_{\xi_j}\in\Gamma(X)$ or $\mu_{\xi_j}\in\Theta(X)$ respectively. 
Therefore we can assume from the beginning that $\widehat\mu_{\xi_j}(y) > 0$.

Since we assume that $\widehat{\mu}_{\xi_j}(y)>0$, it follows from (\ref{Eq1-2}) that

\begin{eqnarray}\label{th1.2.3}
	\widehat{\mu}_{\xi_i}(s,0)=e^{-\sigma_i s^2 }, s\in \mathbb{R},
\end{eqnarray}
where $\sigma_i\geq 0$, for all those $i$ for which $a'_i d'_j - b'_i c'_j\neq 0$ for all  $j=\overline{1,n}$.

Suppose that for some $i$ we have
\begin{eqnarray}\label{th1.2.6.0}
	a'_i d'_j - b'_i c'_j\neq 0 \text{ for all } j=\overline{1,n}
\end{eqnarray}
and at least one of the numbers $a''_i$, $b''_i$ is equal to 1. We can assume without restricting the generality, that    

\begin{eqnarray}\label{th1.2.6}
	a'_1 d'_j - b'_1 c'_j\neq 0 \text{ for all } j=\overline{1,n}
\end{eqnarray}
and $a'_1\neq 0$, $a''_1=1$. 

Substituting $l_1=l$, $l_2=0$ in (\ref{Eq1-1}), we obtain

\begin{eqnarray}\label{Eq1-1.1}
	\prod_{j=1}^{n}  \widehat\mu_{\xi_j} (a'_j s_1+ b'_j s_2, a''_j l )
	=\prod_{j=1}^{n}  \widehat\mu_{\xi_j} (c'_j s_1+ d'_j s_2, c''_j l),\quad s_1,s_2\in\mathbb{R}.
\end{eqnarray}

We can renumber the functions in (\ref{Eq1-1.1}) in such a way that 

\qquad $a'_1 b'_j - b'_1 a'_j\neq 0$, $a''_j=0$ for $j=\overline{2,n_1}$, 

\qquad $a'_1 b'_j - b'_1 a'_j\neq 0$, $a''_j=1$ for $j=\overline{n_1+1,n_2}$, 

\qquad $a'_1 b'_j - b'_1 a'_j=0$, $a''_j=0$ 
for $j=\overline{n_2+1,n_3}$, 

\qquad $a'_1 b'_j - b'_1 a'_j=0$, $a''_j=1$ 
for $j=\overline{n_3+1,n}$.

Suppose that $j>n_3$. Then $\widehat\mu_{\xi_{j}} (a'_j s_1+ b'_j s_2, l )= \widehat\mu_{\xi_{j}} \left({a'_j\over a'_1} (a'_1 s_1+ b'_1 s_2), l  \right)$. Now we can define the characteristic function $$\widehat\mu (s, l )= \widehat\mu_{\xi_{1}} (s, l)  \prod_{j=n_3+1}^{n}  \widehat\mu_{\xi_j} \left({a'_j\over a'_1} s, l  \right).$$

We will find $\widehat\mu (s, 0 )$.
Since $a'_1 b'_j - b'_1 a'_j=0$, we have that $b'_j={b'_1 \over a'_1} a'_j$. If for some $j_0$ the coefficient $a'_{j_0}=0$ then the multiplier $\widehat\mu_{\xi_{j_0}}(s,0)=1$. If for some $j_1$ the coefficient $a'_{j_1}\neq 0$ then $a'_{j_1}d'_j-b'_{j_1}c'_j={a'_{j_1}\over a'_1} (a'_1d'_j-b'_1c'_j)\neq 0$ for all $j=\overline{1,n}$. In this case we have representation (\ref{th1.2.3}). Therefore $\widehat\mu (s, 0 )= e^{-\sigma s^2}$, where $\sigma \geq 0$.

We need some notation. Put

\begin{eqnarray}\label{th1.2.4}
	\psi(s)=\log\widehat{\mu}(s,1),\quad \psi_{a_j}(s)=\log \widehat{\mu}_{\xi_j}(s,a''_j),\quad
	\psi_{c_j}(s)=\log \widehat{\mu}_{\xi_j}(s,c''_j),\quad s\in \mathbb{R}.
\end{eqnarray}

Substituting $l=1$ in (\ref{Eq1-1.1}), we get from (\ref{Eq1-1.1}) that the functions $\psi(s)$, $\psi_{a_j}(s)$ and $\psi_{c_j}(s)$ satisfy the following equation

\begin{eqnarray}\label{th1.2.7}
	\psi(a'_1 s_1+b'_1 s_2)+\sum_{j=2}^{n_3}  \psi_{a_j} \left(a'_j s_1+b'_j s_2 \right)  =
	\sum_{j=1}^n  \psi_{c_j} \left(c'_j s_1+d'_j s_2\right) 
	,\quad s_1,s_2\in \mathbb{R}.
\end{eqnarray}
To solve this equation we will use the finite-difference method. We are going to transform equation (\ref{th1.2.7}) in such a way that the obtained equation will contain only the function $\psi(s)$.

Let $k_n$ be an arbitrary element of $\mathbb{R}$. Substitute $s_1+d'_n k_n$ for $s_1$ and $s_2-c'_n k_n$ for $s_2$ in (\ref{th1.2.7}) and subtract (\ref{th1.2.7}) from the resulting equation. We obtain
\begin{eqnarray}\label{th1.2.8.1}
	\Delta_{(a'_1 d'_n-b'_1 c'_n)k_n} \psi(a'_1 s_1+b'_1 s_2)
	+\sum_{j=2}^{n_3} \Delta_{(a'_j d'_n-b'_j c'_n)k_n} \psi_{a_j} (a'_j s_1+ b'_j s_2) =\nonumber\\=\sum_{j=1}^{n-1} \Delta_{(c'_jd'_n-d'_jc'_n)k_n} \psi_{c_j} (c'_j s_1+ d'_j s_2),\quad s_1,s_2\in \mathbb{R}.
\end{eqnarray}
The right-hand side of equation (\ref{th1.2.8.1}) no longer contains the function $\psi_{c_n}$.

Let $k_{n-1}$ be an arbitrary element of $\mathbb{R}$. Substitute $s_1+d'_{n-1}k_{n-1}$ for $s_1$ and $s_2-c'_{n-1} k_{n-1}$ for $s_2$ in (\ref{th1.2.8.1}) and subtract (\ref{th1.2.8.1}) from the resulting equation. We obtain
\begin{eqnarray}\label{th1.3.1}
	\Delta_{(a'_1d'_{n-1}-b'_1c'_{n-1}) k_{n-1}} \Delta_{(a'_1 d'_n-b'_1 c'_n)k_n} \psi(a'_1 s_1+b'_1 s_2)+ \nonumber\\	
	+\sum_{j=2}^{n_3} \Delta_{(a'_jd'_{n-1}-b'_jc'_{n-1}) k_{n-1}} \Delta_{(a'_j d'_n-b'_j c'_n)k_n} \psi_{a_j} (a'_j s_1+ b'_j s_2) =\nonumber\\=\sum_{j=1}^{n-2} \Delta_{(c'_jd'_{n-1}-d'_jc'_{n-1}) k_{n-1}}
	\Delta_{(c'_jd'_n-d'_jc'_n) k_n} \psi_{c_j} (c'_j s_1+ d'_j s_2),\quad s_1,s_2\in \mathbb{R}.
\end{eqnarray}
The right-hand side of equation (\ref{th1.3.1}) no longer contains the function $\psi_{c_{n-1}}$.

After $n$ steps we obtain the following equation 
\begin{eqnarray}\label{th1.4.1}
	\Delta_{(a'_1d'_1-b'_1c'_1) k_{1}}... \Delta_{(a'_1d'_{n-1}-b'_1c'_{n-1}) k_{n-1}} \Delta_{(a'_1 d'_n-b'_1 c'_n)k_n} \psi(a'_1 s_1+b'_1 s_2)+ \nonumber\\	
	+\sum_{j=2}^{n_3} \Delta_{(a'_jd'_1-b'_jc'_1) k_1} ... 
	\Delta_{(a'_jd'_{n-1}-b'_jc'_{n-1}) k_{n-1}} \Delta_{(a'_j d'_n-b'_j c'_n)k_n} \psi_{a_j} (a'_j s_1+ b'_j s_2) =0 ,\nonumber\\ s_1,s_2\in \mathbb{R},
	\end{eqnarray}
where $k_1, \dots, k_n$ are arbitrary elements of $\mathbb{R}$.

We consider the summands in (\ref{th1.4.1}) where $j=\overline{n_2+1,n_3}$. Since $a'_1 b'_j - b'_1 a'_j=0$, we have that $b'_j={b'_1\over a'_1} a'_j$. If for some $j_0$ the coefficient $a'_{j_0}=0$ then $b'_{j_0}=0$. Therefore the summand $\psi_{a_{j_0}}(s)=0$ and we do not need to remove it from (\ref{th1.4.1}). If for some $j_1$ the coefficient $a'_{j_1}\neq 0$ then $a'_{j_1}d'_j-b'_{j_1}c'_j={a'_{j_1}\over a'_1 }(a'_1d'_j-b'_1c'_j)\neq 0$ for all $j=\overline{1,n}$. In this case we have representation (\ref{th1.2.3}). Since $a''_{j_1}=0$, we have $\psi_{a_{j_1}}(s)=\log\widehat{\mu}_{\xi_{j_1}}(s,0)=-\sigma_{j_1}s^2$. Therefore 
$\Delta_{(a'_{j_1}d'_{1}-b'_{j_1}c'_{1}) k_{1}} ... \Delta_{(a'_{j_1} d'_n-b'_{j_1} c'_n)k_n} \psi_{a_{j_1}} (a'_{j_1} s_1+ b'_{j_1} s_2)$ is a constant. So, we can rewrite (\ref{th1.4.1}) in the following form
\begin{eqnarray}\label{th1.4.2}
	\Delta_{(a'_1d'_1-b'_1c'_1) k_{1}}... \Delta_{(a'_1d'_{n-1}-b'_1c'_{n-1}) k_{n-1}} \Delta_{(a'_1 d'_n-b'_1 c'_n) k_n} \psi(a'_1 s_1+b'_1 s_2)+ \nonumber\\	
	+\sum_{j=2}^{n_2} \Delta_{(a'_jd'_1-b'_jc'_1) k_1} ... 
	\Delta_{(a'_jd'_{n-1}-b'_jc'_{n-1}) k_{n-1}} \Delta_{(a'_j d'_n-b'_j c'_n)k_n} \psi_{a_j} (a'_j s_1+ b'_j s_2)+\nonumber\\  + C =0 ,\quad s_1,s_2\in \mathbb{R},
\end{eqnarray}
where $k_1, \dots, k_n$ are arbitrary elements of $\mathbb{R}$, $C$ is a constant.

Let $l_{n_2}$ be an arbitrary element of $\mathbb{R}$. Substitute $s_1+b'_{n_2}l_{n_2}$ for $s_1$ and $s_2-a'_{n_2}l_{n_2}$ for $s_2$ in (\ref{th1.4.2}) and subtract (\ref{th1.4.2}) from the resulting equation. We obtain
\begin{eqnarray}\label{th1.5.1}
	\Delta_{(a'_1b'_{n_2}-b'_1a'_{n_2}) l_{n_2}} \Delta_{(a'_1d'_1-b'_1c'_1) k_{1}}... \Delta_{(a'_1 d'_n-b'_1 c'_n) k_n} \psi(a'_1 s_1+b'_1 s_2)+ \nonumber\\	
	+\sum_{j=2}^{n_2-1}
	\Delta_{(a'_jb'_{n_2}-b'_ja'_{n_2}) l_{n_2}}
	\Delta_{(a'_jd'_1-b'_jc'_1) k_{1}} ... 
	\Delta_{(a'_j d'_n-b'_j c'_n) k_n} \psi_{a_j} (a'_j s_1+ b'_j s_2)=\nonumber\\  =0 ,\quad s_1,s_2\in \mathbb{R}.
\end{eqnarray}
The left-hand side of equation (\ref{th1.5.1}) no longer contains the function $\psi_{a_{n_2}}$.

After $n+n_2-1$ steps we will obtain the following equation
\begin{eqnarray}\label{th1.6.1}
	\Delta_{(a'_1b'_2-b'_1a'_2) l_2}
	...
	\Delta_{(a'_1b'_{n_2}-b'_1a'_{n_2}) l_{n_2}}
	\Delta_{(a'_1d'_1-b'_1c'_1) k_{1}} ... 
	\Delta_{(a'_1 d'_n-b'_1 c'_n) k_n} \psi (a'_1 s_1+ b'_1 s_2) =0 ,\nonumber\\ s_1,s_2\in \mathbb{R}.
\end{eqnarray}
where $k_1, \dots, k_n$ and $l_2, \dots, l_{n_2}$ are arbitrary elements of $\mathbb{R}$. Since condition (\ref{th1.2.6}) is fulfilled and $a'_1b'_j-b'_1a'_j\neq 0$ for $j=\overline{2,n_2}$,  it follows from  (\ref{th1.6.1}) that the function $\psi$ satisfies the following equation  
\begin{eqnarray}\label{th1.2.9}
	\Delta^{n+n_2-1}_{h} \psi (s)=0,\quad s, h\in \mathbb{R}.
\end{eqnarray}
We obtain from (\ref{th1.2.9}) that
\begin{eqnarray}\label{th1.2.10}
	\widehat\mu(s,1)=e^{-P(s)},
\end{eqnarray}
where $P(s)$ is a polynomial of degree $n$. Taking into account the representation $\widehat\mu (s, 0 )= e^{-\sigma s^2}$, it follows from Lemma \ref{EntireFunction} and its inequality (\ref{EntireFunctionIneq}) that the degree of $P(s)$ is at most 2. Taking into account that $\widehat\mu(-s,1)=\widehat\mu(s,1)$, we obtain
$\widehat\mu(s,1)=\kappa e^{-\sigma' s^2}$. Since $\mu\in M^1(X)$, it follows from Lemma \ref{ClassTheta} that $\mu\in \Theta(X)$. Lemma \ref{AlmostKramer} implies that $\mu_{\xi_1}\in \Theta(X)$.
Similarly, we obtain that $\mu_{\xi_i}\in\Theta(X)$  
for all those $i$ for which $a'_i d'_j - b'_i c'_j\neq 0$ for all  $j=\overline{1,n}$ and at least one of the numbers $a''_i$, $b''_i$ is equal to 1. Thus, we proved statement (2) of Theorem \ref{main.nonzero}.

Now we prove statement (1) of Theorem \ref{main.nonzero}. Note that if $a\in Aut(\mathbb{R}\times \mathbb{Z}(2))$ then $a'\neq 0$ and $a''=1$. 

Suppose that for some $i$ we have $a_i d_j - b_i c_j\in Aut(X)$ for all  $j=\overline{1,n}$. Hence $a_i d_j - b_i c_j\in Aut(Y)$ for all  $j=\overline{1,n}$. We can assume without restricting the generality, that

\begin{eqnarray}\label{th1.0}
	a_1 d_j - b_1 c_j\in Aut(Y) \text{ for all } j=\overline{1,n}.
\end{eqnarray}

It follows from (\ref{th1.0}) that (\ref{th1.2.6}) is fullfilled and at least one of the numbers $a''_1$, $b''_1$ is equal to 1. Hence $\mu_{\xi_1}\in \Theta(X)$ and $\widehat{\mu}_{\xi_{{1}}}$ has form (\ref{def1}).

Substituting $s_1=s_2=0$ in (\ref{Eq1-1}), we obtain
\begin{eqnarray}\label{Eq1-2.1}
	\prod_{j=1}^{n}  \widehat\mu_{\xi_j} (0,a''_j l_1+ b''_j l_2)=\prod_{j=1}^{n}  \widehat\mu_{\xi_j} (0,c''_j l_1+ d''_j l_2 ),\quad l_1,l_2\in\mathbb{Z}(2).
\end{eqnarray}
It follows now from (\ref{th1.0}), Lemma \ref{Equation} and Proposition \ref{z2} that 
\begin{eqnarray}\label{Eq1-2.1}
	\widehat\mu_{\xi_1} (0, l)= 1 ,\quad l\in\mathbb{Z}(2).
\end{eqnarray}
Since $\widehat{\mu}_{\xi_{{1}}}$ has form (\ref{def1}), we obtain now that
$\kappa=1$ in (\ref{def1}). Hence $\mu_{\xi_1}\in \Gamma(X)$. Thus, we proved statement (1) of Theorem \ref{main.nonzero}. 

\hfill $\square$

\begin{remark}\label{cor1H}
	{\rm Let $L_1$, $L_2$ be as in Theorem \ref{main.nonzero}. Suppose that $L_3=L_1$ and $L_4=-L_2$. If the vectors $(L_1,L_2)$ and $(L_1,L_2)$ are identically distributed, then the conditional distribution of $L_2$ given $L_1$ is symmetric. Hence Theorem \ref{main.nonzero}		 implies the following analogue of the Heyde Theorem:}
		
		\textit{Let $X=\mathbb{R}\times \mathbb{Z}(2)$.	
		Let $\xi_1, \dots, \xi_n$ be independent random variables with values in $X$ and distributions $\mu_{\xi_j}$ with non-vanishing characteristic functions. Consider the linear forms $(\ref{th1.1.1})-(\ref{th1.1.2})$, 
		where $a_j, b_j\in End(X)$. Suppose that the conditional distribution of $L_2$ given $L_1$ is symmetric. We have the following possibilities:}
		
		\textit{$(1)$ for all those $i$ for which $a_i b_j + b_i a_j\in Aut(X)$ for all  $j=\overline{1,n}$ distributions $\mu_{\xi_i}\in\Gamma(X)$;}
		
		\textit{$(2)$ for all those $i$ for which $a'_i b'_j + b'_i a'_j\neq 0$ for all  $j=\overline{1,n}$ and at least one of the numbers $a''_i$, $b''_i$ is equal to 1 distributions $\mu_{\xi_i}\in\Theta(X)$;}
		
		\textit{$(3)$ for all those $i$ for which $a'_i b'_j + b'_i a'_j\neq 0$ for all  $j=\overline{1,n}$ distributions $\mu_{\xi_i}\in\Lambda(X)$.}
		
	\rm{Note that if $a_j, b_j \in Aut(X)$ then $a_i b_j + b_i a_j\not\in Aut(X)$ and $a''_i=b''_i=1$ for all $i=\overline{1,n}$. Therefore		
		  statements (1) and (3) are impossible. Statement (2) can be compared with \cite[Theorem 2.1]{Fe2022}.}
	\end{remark}

\begin{remark}\label{cor2SD}
	{\rm We consider two sets of independent random variables $\xi_1, \dots, \xi_n$ and $\xi'_1, \dots, \xi'_n$, where $\xi_j$ and $\xi'_j$ are identically distributed.  Let $L_1$, $L_2$ be as in Theorem \ref{main.nonzero}. We also suppose that $a_j, b_j \in Aut(X)$. Suppose that $L_3=L_1$ and $L_4=L'_2=b_1\xi'_1+\cdots+b_n\xi'_n$. 
	In this case $a_i b_j + b_i a_j\in Aut(X)$ for all $i,j=\overline{1,n}$. therefore statement (1) of Theorem \ref{main.nonzero} is possible because $a_ib_j\in Aut(X)$.	
	Thus, we have the following statement:
		all distributions $\mu_{\xi_i}\in\Gamma(X)$.
		
	If the vectors $(L_1,L_2)$ and $(L_1,L'_2)$ are identically distributed, then $L_1$ and $L_2$ are independent. Therefore the above statement can be considered as an analogue of the Darmois-Skitovich theorem for the group $X=\mathbb{R}\times\mathbb{Z}(2)$ (compare with \cite[\S10.3]{FeBook2}).}
\end{remark}

Consider now the case where the characteristic functions of the distributions can have zeroes.
The following lemma was proved in \cite{My2024}.

\begin{lemma}(\cite{My2024})\label{KlebanovZ2} 
	Let $X=\mathbb{Z}(2)$. Then there exist independent random variables $\xi_j, j=\overline{1,n}, n \ge 2,$ with values in $X$ and distributions $\mu_{\xi_j}$, and endomorphisms $a_j, b_j, c_j, d_j $ such that the random vectors $(L_1,L_2)$ and $(L_3,L_4)$ are identically distributed, but for all those $i$, for which condition $a_i d_j - b_i c_j = 1$ is valid for all $j=\overline{1,n}$, $\mu_i\not\in I(X)$. 
\end{lemma}

Lemma \ref{KlebanovZ2} shows that we can not expect to obtain a reasonable statement in the case $X=\mathbb{R}\times\mathbb{Z}(2)$. Nevertheless, at the beginning of the proof of Theorem \ref{main.nonzero} we obtained representation (\ref{Eq1-2}) without any restrictions on the characteristic functions. Therefore we can formulate the following theorem.  

\begin{theorem}\label{main.zero}
	Let $X=\mathbb{R}\times \mathbb{Z}(2)$.	
	Let $\xi_1, \dots, \xi_n$ be independent random variables with values in $X$. Consider the linear forms $(\ref{th1.1.1})-(\ref{th1.1.4})$, 
	where $a_j, b_j, c_j, d_j\in End(X)$. Suppose that the random vectors $(L_1,L_2)$ and $(L_3,L_4)$ are identically distributed. 
	Then 
	for all those $i$ for which $a'_i d'_j - b'_i c'_j\neq 0$ for all  $j=\overline{1,n}$ distributions $\mu_{\xi_i}\in\Lambda(X)$.
\end{theorem}

\section{Funding}

This research was supported in part by the Akhiezer Foundation.

\section{Acknowledgements}

The author is grateful to Alexander Il’inskii for valuable remarks.



\end{document}